\documentclass{article}
\usepackage{amsmath,amsthm,amsfonts,amscd}
\DeclareMathOperator{\Hom}{Hom}
\DeclareMathOperator{\Com}{Com}
\DeclareMathOperator{\II}{I}
\DeclareMathOperator{\1}{id}
\DeclareMathOperator{\ad}{ad}
\DeclareMathOperator{\Ker}{Ker}
\DeclareMathOperator{\IM}{Im}
\DeclareMathOperator{\ch}{char}
\newcommand{\NN}{\mathbb{N}}
\newcommand{\EEnd}{\mathcal Com}
\newcommand{\EE}{\mathcal E}
\newcommand{\bul}{\bullet}

\newcommand{\w}{\omega}
\newcommand{\p}{\partial}
\renewcommand{\=}{:=}
\renewcommand{\t}{\otimes}
\renewcommand{\o}{\circ}

\newtheorem{thm}{Theorem}
\theoremstyle{definition}
 \newtheorem{defn}[thm]{Definition}
\theoremstyle{definition}
 \newtheorem{exam}[thm]{Example}
\begin{document}
\title{Operadic deformations as a tool for cogravity}
\author{Eugen Paal\\ \\
Department of Mathematics, Tallinn Technical University\\
Ehitajate tee 5, 19086 Tallinn, Estonia\\
e-mail: eugen@edu.ttu.ee
}
\date{}
%
\maketitle
\thispagestyle{empty}
\begin{abstract}
The deformation equation and its integrability condition (Bianchi identity)
of a non-(co)associative deformation in operad algebra are found.
Based on physical analogies, \emph{cogravity} equations are proposed.
\end{abstract}

\section{Introduction}

Non-associativity is sometimes said to be an \emph{algebraic} equivalent of the
differential geometric concept of curvature (e.~g. \cite{NeSa1,NeSa2}). By
adjusting this for physics, one may surmise that gravity and gauge fields
can be described in terms of non-associative algebras. In particular, instead of the
\emph{curvature} of the space-time, \emph{associator} rises to the fore. In
this sense, gravity can be seen to have an algebraic representation. When
non-associativity of space-time becomes large, operadic structure
will become important and one must use \emph{operad algebra} to understand
the algebraic underground of the gravity and how gravity could be
quantized. Instead of the \emph{quantum gravity}, the
\emph{operadic gravity} rises to the fore.

In this paper, the equivalence is clarified from the \emph{linear} deformation
theoretical point of view. By using the Gerstenhaber brackets and a coboundary
operator in a pre-operad, the (formal) associator can be represented
as a curvature form in differential geometry. This (structure) equation is called a
\emph{deformation equation}.
Its integrability condition is the Bianchi identity.

Based on physical analogies, \emph{cogravity} equations are proposed.

\section{Operad algebra}

Let $K$ be a unital associative commutative ring, $\ch K\neq2,3$, and let $C^n$
($n\in\NN$) be unital $K$-modules. For \emph{homogeneous} $f\in C^n$,
$n$ is called the \emph{degree} of $f$ and (when it does not
cause confusion) $f$ is written instead of $\deg f$. For example, $(-1)^f\=(-1)^n$,
$C^f\=C^n$ and $\o_f\=\o_n$. Also, it is convenient to use the
\emph{reduced} degree $|f|\=n-1$. Throughout this paper, it is assumed that
$\t\=\t_K$.

\begin{defn}
A linear \emph{pre-operad} (\emph{composition system}) with coefficients in
$K$ is a sequence $C\=\{C^n\}_{n\in\NN}$ of unital $K$-modules (an
$\NN$-graded $K$-module), such that the following conditions hold.
\begin{enumerate}
\item[(1)]
For $0\leq i\leq m-1$ there exist \emph{(partial) compositions}
\[
  \o_i\in\Hom(C^m\t C^n,C^{m+n-1}),\qquad |\o_i|=0.
\]
\item[(2)]
For all $h\t f\t g\in C^h\t C^f\t C^g$, the \emph{composition
(associativity) relations} hold,
\[
(h\o_i f)\o_j g=
\begin{cases}
    (-1)^{|f||g|} (h\o_j g)\o_{i+|g|}f
                       &\text{if $0\leq j\leq i-1$},\\
    h\o_i(f\o_{j-i}g)  &\text{if $i\leq j\leq i+|f|$},\\
    (-1)^{|f||g|}(h\o_{j-|f|}g)\o_i f
                       &\text{if $i+f\leq j\leq|h|+|f|$}.
\end{cases}
\]
\item[(3)]
There exists a unit $\II\in C^1$ such that
\[
\II\o_0 f=f=f\o_i \II,\qquad 0\leq i\leq |f|.
\]
\end{enumerate}
\end{defn}

In the 2nd item, the \emph{first} and \emph{third} parts of the defining
relations turn out to be equivalent.

Elements of an operad may be called \emph{operations}.
Operad can be seen as a system of operations closed with respect to
compositions.

\begin{exam}[composition pre-operads]
\label{CG}
Let $L$ be a unital $K$-module and
$\EE_L^n\={\EEnd}_L^n\=\Hom(L,L^{\t n})$. Define the partial compositions
for $f\t g\in\EE_L^f\t\EE_L^g$ as
\[
f\o_i g\=(-1)^{i|g|}(\1_L^{\t i}\t g\t\1_L^{\t(|f|-i)})\o f,
         \qquad 0\leq i\leq|f|.
\]
Then $\EE_L\=\{\EE_L^n\}_{n\in\NN}$ is a pre-operad (with the unit
$\1_L\in\EE_L^1$) called the \emph{Cartier (composition) pre-operad} of $L$.
\end{exam}

Thus, \emph{algebraic co-operations} can be seen to be elements of a
composition operad. As \emph{opposed} to this, \emph{algebraic operations}
can be seen to be elements of the \emph{Gerstenhaber (composition) pre-operad}
$\{\Hom(L^{\t n},L)\}_{n\in\NN}$ \cite{Ger}.

 \section{Gerstenhaber brackets}

The \emph{total composition} $\bul\:C^f\t C^g\to C^{f+|g|}$ is defined by
\[
f\bul g\=\sum_{i=0}^{|f|}f\o_i g\in C^{f+|g|},
     \qquad |\bul|=0.
\]
The pair $\Com C\=\{C,\bul\}$ is called the \emph{composition algebra} of $C$.

The \emph{Gerstenhaber brackets} $[\cdot,\cdot]$ are defined in $\Com C$ by
\[
[f,g]\=f\bul g-(-1)^{|f||g|}g\bul f=-(-1)^{|f||g|}[g,f],\qquad|[\cdot,\cdot]|=0.
\]
The \emph{commutator algebra} of $\Com C$ is denoted as
$\Com^{-}\!C\=\{C,[\cdot,\cdot]\}$. It turns out that $\Com^-\!C$ is a
\emph{graded Lie algebra}. The Jacobi
identity reads
\[
(-1)^{|f||h|}[[f,g],h]+(-1)^{|g||f|}[[g,h],f]+(-1)^{|h||g|}[[h,f],g]=0.
\]

In a pre-operad $C$, define a \emph{pre-coboundary} operator $\p_\Delta$ by
\begin{align*}
\p_\Delta f&\=\ad_\Delta^{right}f\=[f,\Delta],\qquad |\p_\Delta|=|\Delta|.
\end{align*}
It follows from the Jacobi identity the (right) derivation property
\[
\p_\Delta[f,g]=(-1)^{|\Delta||g|}[\p_\Delta f,g]+[f,\p_\Delta g]
\]
and the commutation relation
\[
[\p_f,\p_g]:=\p_f\p_g-(-1)^{|f||g|}\p_g\p_f=\p_{[g,f]}.
\]
Thus, if $|\Delta|$ is \emph{odd}, then
\[
\p_\Delta^{2}=\frac{1}{2}[\p_\Delta,\p_\Delta]=
\frac{1}{2}\p_{[\Delta,\Delta]}=\p_{\Delta\bul\Delta}=\p_{\Delta^{2}}.
\]

\section{Deformation equation}

For an operad $C$, let $\Delta,\Delta_0\in C^{2}$.
The difference $\w:=\Delta-\Delta_0$
is called a \emph{deformation}, and $\Delta$ is said to be a deformation of
$\Delta_0$. Let $\p:=\p_{\Delta_0}$, and denote the (formal) associators of
$\Delta$ and $\Delta_0$ as follows:
\[
A:=\Delta\bul\Delta=\frac{1}{2}[\Delta,\Delta],
\qquad
A_0:=\Delta_0\bul\Delta_0=\frac{1}{2}[\Delta_0,\Delta_0].
\]
For a Cartier pre-operad (Example \ref{CG}), the formal associator
reads  as a \emph{coassociator}
\[
A=(\Delta\t\1_L-\1_L\t\Delta)\o\Delta.
\]
The deformation is called \emph{associative} if $A=0=A_0$.

To find the deformation equation, calculate
\begin{align*}
A
&=\frac{1}{2}[\Delta_0+\w,\Delta_0+\w]\\
&=\frac{1}{2}[\Delta_0,\Delta_0]+\frac{1}{2}[\Delta_0,\w]
     +\frac{1}{2}[\w,\Delta_0]+\frac{1}{2}[\w,\w]\\
&=A_0-\frac{1}{2}(-1)^{|\Delta_0||\w|}[\w,\Delta_0]
     +\frac{1}{2}[\w,\Delta_0]+\frac{1}{2}[\w,\w]\\
&=A_0+[\w,\Delta_0]+\frac{1}{2}[\w,\w].
\end{align*}
So we get the \emph{deformation equation}
$$
\boxed{A-A_0=\p\w+\frac{1}{2}[\w,\w]}
$$
The deformation equation can be seen as a differential equation for
$\omega$ with given associators $A_0, A$. Note that if the associator is
fixed, i.~e. $A=A_0$, one obtains the \emph{Maurer-Cartan (master)
equations}, well-known from the theory of associative deformations.

\section{Prolongation}

Now differentiate the deformation equation,
\begin{align*}
\p(A-A_0)
&=\p^{2}\w+\frac{1}{2}\p[\w,\w]\\
&=\p^{2}\w+\frac{1}{2}(-1)^{|\p||\w|}[\p\w,\w]
     +\frac{1}{2}[\w,\p\w]\\
&=\p^{2}\w-\frac{1}{2}[\p\w,\w]
     +\frac{1}{2}[\w,\p\w]\\
&=\p^{2}\w-\frac{1}{2}[\p\w,\w]
     -\frac{1}{2}(-1)^{|\p\w||\w|}[\p\w,\w]\\
&=\p^{2}\w-[\p\w,\w].
\end{align*}
Again using the deformation equation, we obtain
\begin{align*}
\p(A-A_0)
&=\p^{2}\w-[\p\w,\w]\\
&=\p^{2}\w-[A-A_0-\frac{1}{2}[\w,\w],\w]\\
&=\p^{2}\w-[A-A_0,\w]+\frac{1}{2}[[\w,\w],\w].
\end{align*}
Finally use $[[\w,\w],\w]=0$ to obtain the condition
$$
\boxed{\p(A-A_0)=\p^{2}\w-[A-A_0,\w] }
$$

\section{Associativity constraint and Bianchi identity}

We know that $\p^{2}=\p_{A_0}$. Hence, if \emph{associativity} constraint $A_0=0$ holds,
then
$$
\boxed{\p^{2}=0}
$$
In the case of the Cartier pre-operad,
the coset $\Ker\p/\IM\p$ is called
the \emph{Cartier cohomology} \cite{Ca,SnSt}.
The deformation equation for such a \emph{non-associative} deformation
reads
$$
\boxed{A=\p\w+\frac{1}{2}[\w,\w] }
$$
One can see that associator is a formal \emph{curvature} while the
deformation is working as a \emph{connection}. One can say that associator
is an \emph{operadic} equivalent of the curvature.
The integrability condition of the deformation equation reads as the
\emph{Bianchi identity}
\[
\boxed{\p A=[\w,A]}
\]
One can easily check that further differentiation does not add new conditions.

\section{Covariant derivation}

Note that
\[
\p_\Delta f
=[f,\Delta]
=[f,\Delta_0+\omega]
=[f,\Delta_0]+[f,\omega]=\p f+[f,\omega].
\]
One can say that $\nabla:=\p_\Delta$ is a \emph{covariant derivation}.
The Bianci identity reads
\[
\boxed{\nabla A=\p A+[A,\omega]=0}
\]
Also note that
\[
\boxed{\nabla^{2}f=[f,A]}
\]
So the condition $\nabla^{2}=0$ does not imply $A=0$. Instead,
$\nabla^{2}=0$ implies that $A$ lies in the \emph{center} of
$\Com^{-}C$. In particular,
\[
\nabla^{2}=0\quad\Longrightarrow
\quad\p A=0
\quad\Longrightarrow\quad
A\in\Ker\p.
\]
Note that $\Delta$ may nevertheless remain \emph{non-(co)associative}.

\section{Cogravity versus gravity}

Thus the differential geometrical notion of curvature can be easily
adjusted for deformations in a pre-operad. Rather than to speak about
\emph{algebraic} deformation theory, one may speak about the \emph{geometrical}
one.
Geometry performs the pioneering role in creating of
the exact scientific world picture.
One may ask that how far one can proceed with geometrical notions
in deformation theory.
In particular, this question may be adjusted for physics as well.

In General Relativity, gravity is a fundamental interaction associated with
the space-time curvature (associator).
As opposed to this, the hypothetic kind of
interaction associated with a coassociator may be called the \emph{cogravity}.
In a sense, cogravity may be called the \emph{opposite gravity} as well.

Let $C$ be a composition operad of a \emph{free} module $L$. For a given
base $e^{i}\in L$, the \emph{connection} (deformation) coefficients
$\Gamma^{i}_{jk}$ and the
coassociator components $A_{jkl}^{i}$ are introduced by
\[
\omega e^{i}=\Gamma^{i}_{jk}e^{j}\t e^{k}, \qquad Ae^{i}=A^{i}_{jkl}e^{j}\t
e^{k}\t e^{l}.
\]
To obtain the \emph{cogravity equations}, one must compose the
\emph{Ricci coassociator} $A_{ij}:=A_{isj}^{s}$ and the
\emph{energy-momentum} tensor $T_{ij}$. Then the  cogravity
equations reads
\[
\boxed{ A_{ij}=k\left(T_{ij}-\frac{1}{2}g_{ij}T\right)}
\]
Dynamical deformations of the \emph{space-time geodesic} \emph{comultiplication}
can be seen as a prospective model to specify (represent) the cogravity equations.
Geodesic multiplication \cite{Ki,Ak} $m$ is a non-associative deformation
of the \emph{vector addition} operation $m_0$. Comultiplications $\Delta,\Delta_0$
are defined for a function $F$ on the space-time by
\[
\Delta F=F\o m,\qquad \Delta_0 F=F\o m_0.
\]

The second idea follows the Maxwell (gauge field) equations.
It is well-known that the \emph{first} pair of the Maxwell equations
can be represented as the Bianchi identity.
To introduce the \emph{second} pair, one must define a \emph{dualization} $^{\dag}$
and an \emph{operad current} $J$. Then the (gauge field) Maxwell-like equations read
\[
\boxed{
\nabla A=\p A+[A,w]=0,\qquad
\nabla A^{\dag}=\p A^{\dag}+[A^{\dag},w]=J^{\dag}
}
\]
This idea seems to be similar to \cite{MaTa}.
A concept of a cogravity was posed in \cite{Ma} as well.
In the above approach, one must study the physical ((co)gravity) equations in
deformation complexes.

\section*{Acknowledgements}

I would like to thank  P.~Kuusk, J.~L\~{o}hmus and J.~Stasheff
for helpful comments.
Research was supported in part by the ESF grant 3654.

\end{document}